\documentclass[12pt]{amsart}

\usepackage{amsfonts,amscd,amsmath,amssymb,amsthm,latexsym,xypic}
\usepackage{graphicx}

\newcommand{\ZZ}{{\bf Z}}
\newcommand{\QQ}{{\bf Q}}
\newcommand{\RR}{{\bf R}}
\newcommand{\CC}{{\bf C}}
\newcommand{\n}{\noindent}
\newcommand{\NN}{{\bf N}} 

\newcommand{\bbP}{\mathbb{P}}
\newcommand{\Pic}{\textup{Pic}}
\def\f{{\varphi}}
\def\l{{\mathcal {L}}}
\def\r{{\Lambda}}
\newtheorem{df}{Definition}[section]
\newtheorem{teo}[df]{Theorem}  

\newtheorem{lemma}[df]{Lemma}

\newtheorem{oss}[df]{Remark}

\begin{document}

\title{Correspondences between $K3$ surfaces}

\author
{Federica Galluzzi} 
\author{ Giuseppe Lombardo}

\thanks{The two authors are supported by EAGER - European Algebraic Geometry Research training network, Contract
number HPRN-CT-2000-00099
and  by Progetto di Ricerca Nazionale COFIN 2000 "Geometria delle Variet\`a
Algebriche"}

\address{Dipartimento di Matematica, Universit\`a di Torino
 Via Carlo Alberto 10 10123 Torino,ITALY}
\email {federica.galluzzi@unito.it}
\address{Dipartimento di Matematica, Universit\`a di Torino
 Via Carlo Alberto 10 10123 Torino,ITALY}
\email{lombardo@dm.unito.it}
\address{Department of Mathematics, University of Michigan, Ann Arbor, MI
48109,USA}
\email{idolga@umich.edu}

\maketitle 

\centerline{with an Appendix by Igor Dolgachev}

\begin{abstract}
In this paper we show that there is a correspondence between some $K3$ surfaces with non-isometric transcendental lattices constructed as a twist of the transcendental lattice of the Jacobian of a generic genus 2 curve.
Moreover, we show the existence of a correspondence between a general $K3$ surface with $\rho =17$ and a Kummer surface having transcendental lattices $\QQ$-Hodge isomorphic.

\end{abstract}

\maketitle

\section{Introduction}
In this paper we study the existence of correspondences between $  K3$ surfaces $ X(k,m,n),\,$ with $k,m,n \in \NN ,\,$  Picard number 17 and transcendental lattices   $T(k,m,n) \cong U(k)\oplus U(m)\oplus \langle -2n \rangle \,$. In the fundamental paper \cite{Mu} Mukai showed that  correspondences between $K3$ surfaces exist if the transcendental lattices are Hodge isometric over $\QQ .\,$ This construction holds if the Picard number of the surfaces is  greater or equal to $ 11 .\,$ Nikulin improved the result afterwards in \cite{N2} obtaining the lower bound $5$ for the Picard number. 

\n The aim of our work is to realize examples of $K3$ surfaces with transcendental lattice which are not Hodge isometric but such that a correspondence between them already exists. This in particular implies the existence of an algebraic cycle 
on the middle cohomology of the product of two surfaces arbitrarily chosen in the constructed family.

\n
More in details, in the first sections we recall some basic notions and results on lattices and correspondences. 
In Section \ref{si} we consider a generic genus 2 curve and we show the existence of a correspondence between the Jacobian of the curve and a $K3$ surface  with isomorphic transcendental lattice. Since this construction involves a second $K3$ surface whose transcendental lattice has quadratic form multiplied by 2, in Section \ref{emb} and \ref{abs} we generalize this first example. First, we construct $K3$ surfaces "twisting" each direct summand of the transcendental lattice of the Jacobian by natural numbers. Then we find correspondences between them using both Mukai's theorem and Shioda-Inose structures which translate the problem into a problem of looking for isogenies between abelian varieties.

\n
In this way we prove in Theorem\ref{nm} that all the $K3$ surfaces $X(k,m,n)$ are in correspondence to each other.

\n
Finally, in Theorem \ref{kummer} we show the existence of a correspondence between a general $K3$ surface of Picard rank 17 and a Kummer surface of the same rank having transcendental lattices $\QQ$-Hodge isomorphic.

\section{Preliminary notions.} \subsection{Definitions.}A {\em lattice} is a free $\ZZ$-module $\l$ of finite rank  with a $\ZZ$-valued symmetric bilinear form $b_{\l}(x,y).\,$ A lattice  is called {\em even} if the quadratic form associated to the bilinear form has only even values, {\em odd} otherwise. A very useful invariant (under base change) of a lattice is its {\em discriminant} $d(\l)$, defined as the determinant of the matrix of its bilinear form. Then, a lattice is called {\em non-degenerate} if the discriminant is non-zero and {\em unimodular} if the discriminant is $\pm 1 .\,$ If the lattice $\l$ is non-degenerate, the pair $(s_+, s_-),\,$ where $s_{\pm}$ denotes the multiplicity of the eigenvalue $\pm 1$ for the quadratic form associated to $\l \otimes \RR ,\,$ is called {\em signature} of $\l .\,$ Finally, we call $\, s_+ +s_- \,$ the {\em rank} of $\l$ and $\, s_+-s_- \,$ its {\em index}, moreover a lattice is {\em indefinite} if the associate quadratic form has both positive and negative values.

\n 
Given a lattice $\l \,$ we can construct the lattice $\l (m) ,\,$ that is the $\ZZ$-module  $\l $ with bilinear form $b_{\l (m)}(x,y)=mb_{\l}(x,y).$  

\n
An {\em isometry} of lattices is an isomorphism preserving the bilinear form. 
Given a sublattice $\l \hookrightarrow \l ^{\prime},\,$ the embedding is {\em primitive} if $\displaystyle {\l ^{\prime} \over \l }$
is free.

\subsection{Examples.} 
\begin{enumerate}
	\item [i)]
The lattice $\langle n \rangle $ is a free $\ZZ$-module of rank one , $\ZZ  \langle e \rangle ,\,$ with bilinear form \break $b(e,e)=n$.

\item [ii)] The {\em hyperbolic lattice}  is the even, unimodular, indefinite lattice with $\ZZ$-module $\ZZ \langle  e_1,e_2 \rangle $ and bilinear associated form of matrix $\left( \begin{matrix}0&1 \cr 1&0 \end{matrix}\right)$.

\item [iii)]The lattice $E_8$ has  $\ZZ^8$ as  $\ZZ$-module and the matrix of the bilinear form is the  Cartan matrix of the root system of $E_8$. It is an even, unimodular and positive definite lattice.

\end{enumerate}

\subsection{$K3$ and tori lattices.} If $X$ is  a $K3$ surface, one can show that $H^2(X, \ZZ)$ is free of rank $22$ and that there is an isometry $H^2(X, \ZZ)\cong U^3\oplus (E_8(-1))^2$. From now on we denote with $\r$ this $K3$-lattice. For $X$ a complex torus, one has  $H^2(X, \ZZ)\cong U^3 .$

\subsection{Hodge structures.} 
Let $X$ be an abelian or $K3$ surface. If we consider the Hodge decomposition of $H^2(X, \CC)=H^{2,0}(X)\oplus H^{1,1}(X)\oplus H^{0,2}(X)$, inside $H^2(X,\ZZ)$ one has two sublattices, the {\em N\'eron-Severi lattice}
$$
NS(X)\, := \,H^2(X,\ZZ)\cap H^{1,1}(X)
$$ 
and the orthogonal complement of  $NS(X),\,$ the {\em transcendental lattice}  $T_X$ which   has a natural Hodge structure induced by the one of $H^2(X,\ZZ )$.
The {\em  Picard number} of $X,\,$ denoted by $\rho (X),\,$ is the rank of $NS(X).$

\n
A {\em Hodge isometry} between the transcendental lattices of two $K3$ (or abelian) surfaces is an isometry preserving the Hodge decomposition.

\subsection{Fourier-Mukai partners.}
Let $X,Y$ be both abelian (or $K3$) surfaces, they are called {\em Fourier-Mukai partners }   if their derived categories of bounded complexes of \break coherent sheaves are equivalent. Mukai and Orlov have showed that 
\begin{teo}\cite{Mu,Or} $X$ and $Y$ as above are Fourier-Mukai partners if and only if their transcendental lattices are Hodge-isometric.
\end{teo}

\section{Known results.} 
In this section we recall some fundamental results which will play a key role in the next paragraphs. First of all, since for a $K3$ or an abelian surface the period map is surjective,  it is possible to prove the following

\begin{teo}\cite[Corollary 1.9]{M}\label{sur}
Suppose $T \hookrightarrow U^3$ (resp. $T \hookrightarrow \r$) is a primitive sublattice of signature $(2,4-\rho )$ (resp. $(2,20-\rho )$). Then there exists an abelian surface (resp. algebraic $K3$ surface) $X$ and an isometry $T_X \stackrel{\sim}{\rightarrow} T$.

\end{teo}

\n
\begin{df} {\em A $K3$ surface $X$ admits a {\em Shioda-Inose structure} if there is an involution $\iota$ on $X$ such that $\iota ^* (\omega)=\omega$ for every $\omega \in H^{2,0}(X)$ ($\iota$ is called a {\em Nikulin involution}) and with rational quotient map $\pi : X --> Y $ where $Y$ is a Kummer surface. The map $\pi _* $ induces a Hodge isometry $T_{X}(2) \cong T_Y .\,$ This gives a diagram}
\end{df}

\smallskip
\centerline{
\xymatrix{
X \ar@{-->}[rd] & & Z \ar@{-->}[ld] \\
& Y &
}}

\n
of rational maps of degree $2 \,$ where $Z$ is a complex torus and $Y$ is the Kummer surface of $Z.\,$

\n One can detect the existence of a Shioda-Inose structure on a $K3$ surface analyzing the transcendental lattice of the surface.

\begin{teo}\cite[Theorem 6.3]{M}\label{iff} Let $X$ be an algebraic $K3$ surface. $X$ admits a Shioda-Inose structure if and only if there is a primitive embedding $T_X \hookrightarrow U^3 .\,$
\end{teo}

\begin{oss}
{\em Obviously, the Shioda-Inose structure realizes a correspondence between the $K3$ surface $X$ and the abelian variety $Z$; this correspondence will play a fundamental role in 
the following paragraphs.}

\end{oss}

\section{Starting problem.}\label{si}

Let $C$ be a generic genus 2 curve (i.e. such that its Jacobian surface has $\rho (JC)=1$); $JC$ is a principally polarized abelian variety and, if $E$ is the principal polarization, so we have $E^2=2$ then $T:=T_{JC}=U^2\oplus \langle -2 \rangle$. Since we have obvious primitive embeddings of $T$ in $\Lambda$, from Theorem \ref{sur} there exists an algebraic $K3$ surface $X_1$ such that  $T_{X_1} \cong T_{JC}$. Moreover, by \cite[Proposition 6.2]{Mu}, the number of Fourier-Mukai partners of $X_1$ is one since $\rho(X_1)=17$, so such a $X_1$ is  unique (up to isomorphisms). We analyze now the relations between the Jacobian surface and this $K3$ surface. 
First, we observe that we can construct an embedding $T_{X_1}\hookrightarrow U^3$ in the following way; we send the first two copies of $U\subset T_{X_1}$ to the corresponding ones of $U^3$ and the element of square $-2$ to $e_1^3-e_2^3$ where $\{e_j^i\}^{i=1,2,3}_{j=1,2},\,$ is a basis of $U^3$. 

\n Then, from Theorem \ref{iff} $X_1$ admits a Shioda-Inose structure

\bigskip

\centerline{
\xymatrix{
X_1 \ar@{-->}[rd] & & Z \ar@{-->}[ld] \\
& X_2  &
}}

\bigskip

\n and one has $T_Z \cong T_{X_1}$, thus $T_Z\cong T_{JC}.\,$
Recently, Hosono, Lian, Oguiso and Yau have proved the following

\begin{teo} \cite[Main Theorem, 1)]{H-L-O-Y} \label{hloy} 
Let $A$ and $B$ abelian surfaces, they are Fourier-Mukai partners if and only if $Kum(A) \cong Kum(B)$.
\end{teo}

\n It follows that $Kum(Z)\cong Kum(JC)$ since $Z$ and $JC$ are Fourier-Mukai partners. In this way we obtain a correspondence between $X_1,\, X_2$ and $JC$

\bigskip

\centerline{
\xymatrix{
JC \ar@{-->}[rd] &  & Z\ar@{-->}[rd] & &X_1\ar@{-->}[ld] \\
& Kum(JC) & \cong & X_2=Kum(Z) &  
}}
\begin{oss}
{\em We observe that the Hodge isometry between $T_{JC}$ and $T_Z$ can be extended to a Hodge isometry between the second cohomology groups since $JC$ and $Z$ are principally polarized abelian surfaces of rank one. From \cite{Sh}, we obtain that $Z\cong JC$ or $Z\cong {(JC)}^{\vee}$.}
\end{oss}
\section{Embeddings of twisted lattices.}\label{emb}

We want to generalize the situation of the previous paragraph.
We consider the lattice $T(k,m,n):=U(k) \oplus U(m) \oplus \langle -2n \rangle \,,$ with $ k,m,n \in \NN ,\, $ obtained by twisting the summands of $T\,$, equipped with the Hodge structure induced by $T$ itself. 
This is an indefinite even lattice of signature $(2,3).\,$
The aim is to prove the existence of correspondences between $K3$ surfaces having transcendental lattice Hodge-isometric to $T(k,m,n).\,$      
In order to construct such surfaces, we start proving the following

\smallskip

\begin{lemma}
The lattice $T(k,m,n)$ is a primitive sublattice of the $K3$-lattice  $\Lambda .\,$
\end{lemma}

\begin{proof}
First we observe that a lattice $\langle -2t \rangle ,\, t \in \NN ,\, $ can be primitively embedded in $E_8(-1)$ from \cite[Theorem 1.12.4]{N1}. We denote by $\theta _t \in E_8(-1)$ its generator. 

\n In order to construct the required primitive embedding of $T(k,m,n)$ in $\Lambda$ we consider  a $\theta _k$ contained in the first copy of $E_8(-1)$ and a $\theta _n$ in the second one.

\n Then, if  $\{ \delta_1, \dots , \delta_5\}$ is the standard basis  of  $T(k,m,n)$,  we can obtain the desired embedding in the following way:
$$\begin{array}{ccl}
T(k,m,n)& \hookrightarrow &  U^3\oplus E_8(-1)^2 \\
 & & \\
\delta_1 & \mapsto & (1, k)+( 0,0)+(0,0)+\theta_k \\
\delta_2 & \mapsto & (1, 0)+( 0,0)+(0,0) \\
\delta_3 & \mapsto &(0,0)+(1,0)+(0,m) \\
\delta_4 & \mapsto &(0,0)+( 0,0)+(1,0) \\
\delta_5 & \mapsto & \theta_n.
\end{array}$$
\end{proof}

\n
\begin{oss}{\em 
This Lemma, together with Theorem \ref{sur}, implies the existence, for any \break $k,m,n \in \NN$, of a $K3$ surface $X(k,m,n)$  with a Hodge isometry between the transcendental lattices $T_{X(k,m,n)}$  and $T(k,m,n)\,$.  Such a surface is unique (up to isomorphisms) by \cite[Proposition 6.2]{Mu} .}
\end{oss} 

\smallskip
\n
In order to generalize the example constructed in the previous paragraph, we want to produce correspondences among these $K3$ surfaces. We observe that, with these notations, the surfaces of the previous section can be rewritten as $X_1= X(1,1,1)$ and $X_2=X(2,2,2)$. 

\medskip
\n
The problem of finding correspondences among $K3$ surfaces was investigated by Mukai in \cite{Mu}. His idea was to construct such correspondences starting from Hodge isometries (over $\QQ$) between the transcendental lattices of the surfaces. This method works if  the Picard number of the surfaces is sufficiently large. Mukai's result is the following

\n
\begin{teo}\cite[Cor.1.10.]{Mu}\label{tmuk}
Let $X,Y$ be $K3$ surfaces with $\rho (X),\rho(Y) \geq 11 .\,$ If  \break
$\f:T_X\otimes \QQ \rightarrow T_Y\otimes \QQ$ is a Hodge isometry,  then $\f$ is induced by an algebraic cycle.
\end{teo}

\smallskip

Our $K3$ surfaces, unfortunately, do not satisfy the condition of Mukai's theorem, since an isometry between $T(k,m,n)$ and $T(k',m',n')$ doesn't exist, even over $\QQ .\,$  However, in the same article Mukai proved the Oda's conjecture, as modified by Morrison in \cite{M}, which realizes a correspondence between a $K3$ surface $X$ and an abelian surface provided that the $\QQ$-transcendental lattice $T_X \otimes \QQ$ admits an embedding
in $U^3 \otimes \QQ.\,$ In order to construct the correspondences among the $X(k,m,n) ,\,$ we follow a similar approach, translating the problem in another one involving some (abelian) surfaces which are in correspondence with the given ones.

\section{Abelian surfaces and correspondences.}\label{abs}

Motivated by the work of Morrison and Mukai, we start with the following

\smallskip
\n
\begin{lemma}\label{phi}
There is an embedding of $\QQ -$lattices $\phi : T(k,m,n) \otimes \QQ \hookrightarrow U^3 \otimes \QQ .$
\end{lemma}

\begin{proof}
Let $\{e^i_1,e^i_2 \},\,i=1,2,3 \,$ be the basis of the i-th copy of the hyperbolic lattice $U$ and let $\{ a_1,b_1,a_2,b_2,c \}$ be the basis of $T(k,m,n)$. We define $\phi$ in the following way

\[
\begin{array}{ccl}
a_1 &\longmapsto & e_1^1 \otimes 1 \\
&&\\
b_1 &\longmapsto & e_2^1 \otimes k \\
&&\\
a_2 &\longmapsto & e_1^2 \otimes 1 \\
&&\\
b_2 &\longmapsto & e_2^2 \otimes m \\
&&\\
c &\longmapsto & e_1^3 \otimes 1 + e_2^3 \otimes (-n) . \\
\end{array}
\]
\end{proof}

\smallskip
\n
The existence of such an embedding allows us to prove the following

\begin{teo}\label{k3ab}
For any $k,m,n \in \NN$ there exist abelian surfaces $A_n$ and correspondences
\end{teo}  

\[
 \begin{CD}
  Z(k,m,n) @>>> X(k,m,n)     \\
@VVV  @.    \\
A_n @. 
 \end{CD} 
\]
\begin{proof}
Let us consider the embedding  $\phi :   T(k,m,n) \otimes \QQ \hookrightarrow U^3 \otimes \QQ \,$ of Lemma \ref{phi} and the lattice $T_n := \phi (T(k,m,n)\otimes \QQ) \cap U^3 =U^2\oplus  \langle -2n \rangle \,$ with the Hodge structure induced by $\phi .\,$
Since $T_n$ is a primitive sublattice of $\Lambda ,\, $ there exists a $K3$ surface $Y_n$ and a Hodge isometry $T_{Y_n}\cong T_n .\,$ Consider now the basis of $T_{X(k,m,n)}$  given in Lemma \ref{phi}: the multiplication by $1 \over k$ on the sublattice $\ZZ \langle b_1\rangle$ and by $1 \over m$ on $\ZZ \langle b_2\rangle$ induces a Hodge isometry $T_{X(k,m,n)} \otimes \QQ \cong T_{Y_n} \otimes \QQ \, $ and Theorem \ref{tmuk} gives a correspondence  $Z^{\prime}(k,m,n)$ between $X(k,m,n)$  and $Y_n .\,$  

\n
We can define also a primitive embedding $ T_n \cong T_{Y_n} \hookrightarrow  U^3 \,$ sending $U^2 \subset T_n$ to the first two copies of $U$ in $U^3$ via the identity  and sending the element of square $-2n$ to $e^3_1-ne^3_2 .\,$
By Theorem \ref{iff} the existence of such an embedding is equivalent to the existence of a Shioda-Inose structure on the $K3$ surface  $Y_n .\,$
So there exist an abelian surface $A_n ,\,$  a Hodge isometry $T_{A_n}\cong T_n \,$ and a correspondence $Z^{\prime\prime}_n$ between $Y_n$ and $A_n .\,$ The composition of the two correspondences

\smallskip
\centerline{
\xymatrix{
 &  \ar@{->}[ld]Z^{\prime}(k,m,n)\ar@{->}[rd] & &\ar@{->}[ld] \qquad Z^{\prime \prime}_n \qquad \ar@{->}[rd] \\
 X(k,m,n) & & Y_n& & A_n  
}}

\medskip
\n
gives the desired one. 
\end{proof}


\medskip

\n Now, we are able to prove the following 

\smallskip

\begin{teo}\label{nm}
 Let $X(k,m,n)$, $k,m,n \in\NN$, be a $K3$ surface with transcendental lattice Hodge-isometric to $ T(k,m,n)$. For any $k',m',n'\in \NN$ there exist correspondences

\[
 \begin{CD}
 Z_{k,m,n}^{k',m',n'}  @>>> X(k,m,n)     \\
@VVV @.    \\
 X(k',m',n') @. 
 \end{CD} 
\]
and every $X(k,m,n)$ has a correspondence with the Jacobian surface $JC$ of Section \ref{si}.
\end{teo}
\begin{proof}
Let $Z$ be a principally polarized abelian surface with $\rho(Z)=1$ and let $ \{ \lambda_1, \lambda_2, \lambda_3 , \lambda_4 \}$ be a symplectic basis for the lattice of $Z\,$. Let now, for any $n$, $\, B_n:=\CC^2 / \Lambda_n$ be the complex torus with lattice $\Lambda_n=\ZZ \langle \lambda_1, n\lambda_2, \lambda_3 , \lambda_4\rangle$; the polarization $E$ defines on each $B_n$ a polarization of type $(1,n)$; moreover the abelian surfaces $B_n$ are all obviously isogenous to each other. 

\n On the other hand, one has that $T_{B_n}\cong U^2\oplus \langle -2n \rangle $, so every $B_n$ is a Fourier-Mukai partner of the surface $A_n$ constructed in Theorem \ref{k3ab}. This means, by Theorem \ref{hloy}, that \break $Kum(A_n)\cong Kum(B_n),\,$  thus there is a correspondence between $A_n$ and $B_n .\,$
In this way we have constructed a correspondence between $A_n$ and $A_m ,\,$ for any $n,m$ :

\[
\xymatrix{
A_n \ar@{-->}[d] & &  B_n \ar@{-->}[d]   \ar[r]^{\sim}  & B_m \ar@{-->}[d] & &  A_m \ar@{-->}[d] \\
Kum(A_n) &\cong & Kum(B_n) & Kum(B_m)& \cong & Kum(A_m)
}
\] 

\bigskip
\n
But Theorem \ref{k3ab} gives a correspondence between $X(k,m,n)$ and $A_n$ for any $n$,  so the desired $Z_{k,m,n}^{k',m',n'}$ is the composition of these correspondences. Moreover, there is a correspondence  between $X(1,1,1)$ and $JC$ from Section \ref{si}, thus the statement follows.
\end{proof}

\begin{oss}
{\em We observe that the existence of the correspondence is independent from the chosen Hodge structure on $T$. However, it is necessary not to change the structure when we "twist" the starting lattice $T$ by $(k,m,n)$. This allows us to obtain, for any chosen Hodge structure $H$ on $T$, a family of $K3$ surfaces $X_H(k,m,n)$ in correspondence to each other.}

\end{oss}

Morrison, in \cite[Cor 4.4]{M} showed  that a $K3$ surface $K$ of Picard rank 17 is a Kummer surface if and only if there is an even lattice $T'$ with $T_K\cong U(2)^2\oplus T'(2)$. Then, this allows us to prove the following theorem, which generalizes Theorem 
\ref{nm} 

\begin{teo}\label{kummer} Let $X$ be a general $K3$ surface with $\rho(X)=17$ and  $T_X \otimes \QQ$ isomorphic as a Hodge structure to $T_K\otimes \QQ$, where $K$ is a Kummer surface. Then, there is a correspondence between $X$ and $K$. \end{teo}
\begin{proof}
Let $\psi_K=U^2(2)\oplus\langle -4n \rangle $ be the polarization on $T_K$, we show that there is a $\QQ$-basis of $T_X \otimes \QQ$ in which the polarization is $\psi_X =a \psi_K$ with $a \in \QQ$. Since $\psi_X \in Sym^2(T_X\otimes \QQ)^{MT(T_X\otimes \QQ)}$ (where $MT$ denotes the Mumford-Tate group) it sufficies to show that this space has dimension one. This is a consequence of Schur's Lemma, since $MT(T_x\otimes \QQ)(\CC)=SO(5)(\CC)$ and the action on the Hodge structure is irreducible. So, we have the Hodge $\QQ$-isometries $T_X\otimes \QQ \cong a(U^2(2)\oplus\langle -4n \rangle)\cong U^2\oplus \langle -4na_1a_2\rangle$, where $a={a_1\over a_2}$. From the surjectivity of the period map, the last one is a transcendental lattice of a $K3$ surface $X(1,1,-4na_1a_2)$, which is in correspondence with $X$ from Mukai's work. The statement now follows from Theorem \ref{nm}.\end{proof}

\centerline{\textsc{Acknowledgments}}


\noindent

We would like to thank Prof. Bert van Geemen for having introduced us to this problem and for the useful discussions and valuable suggestions. We also thanks Prof. Igor Dolgachev for his interest in this work.


\appendix

\n
\section {A geometric correspondence between $X_1$ and a Kummer surface}

\n
We use the following result of V. Nikulin \cite{N3}.

\begin{teo} Let $X$ be a K3-surface with $T_X \cong U\oplus \langle -2\rangle$. Then $\Pic(X)$ has only finitely many smooth rational curves which form the following graph.

\begin{figure}[h]
\begin{center}
\includegraphics[width=4in]{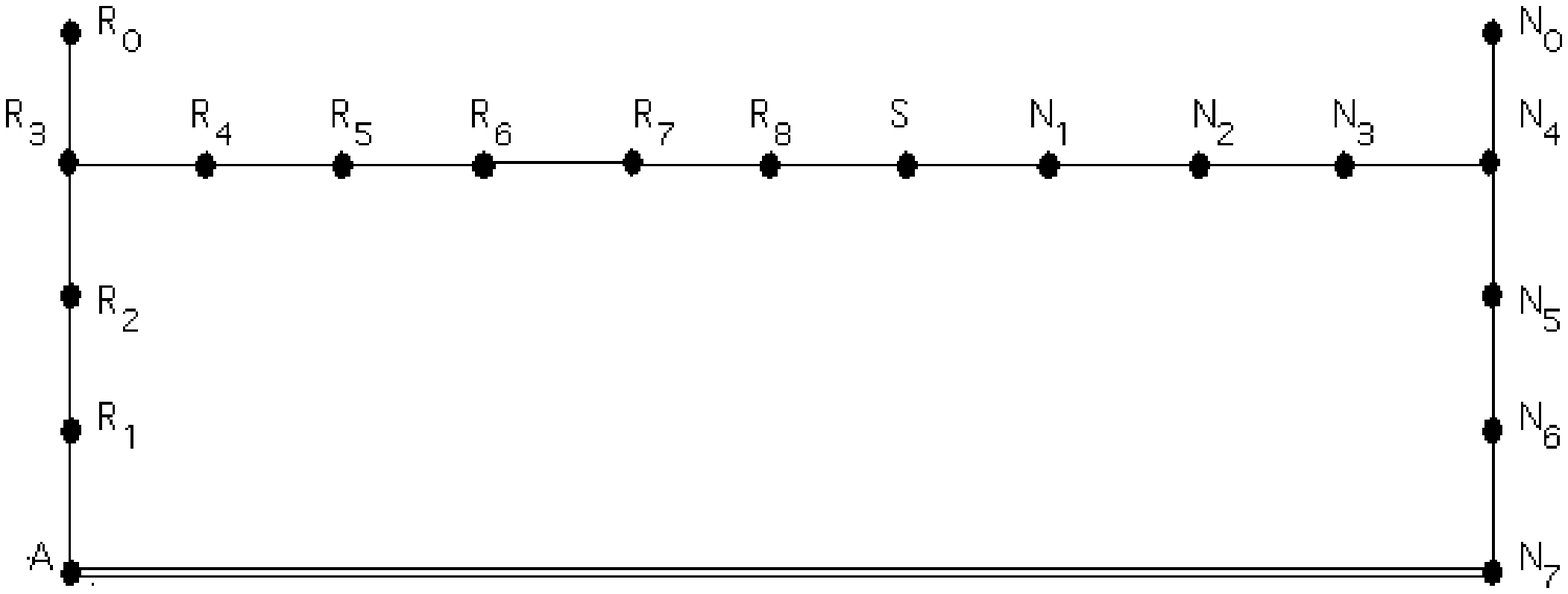}
\label{fanoplane}
\end{center}
\end{figure}

The automorphism group of $X$  is generated by two commuting involutions $\sigma$ and $\tau$. The involution $\sigma$ has 8 isolated fixed points. The set of fixed points of the involution $\tau$ is the union of 8 smooth rational curves and a smooth curve of genus 2. 
\end{teo}

It follows from observing the graph that $X$ admits an elliptic fibration $|F|$ with two singular fibres $$F_1 = 3R_0+2R_1+4R_2+6R_3+5R_4+4R_5+3R_6+2R_7+R_8$$ and 
$$F_2 = 2N_0+N_1+2N_2+3N_3+4N_4+3N_5+2N_6+N_7$$ of  type $\tilde{E}_8$ and $\tilde{E}_7$. It is also has a section $S$. The fixed locus of $\tau$ consists of smooth rational curves $R_1,R_3,R_5,R_7, N_2,N_4,N_6, S$ and a genus 2 curve $W$ which intersects $R_0,N_0, N_7$ with multiplicity 1.  Let $p = W\cap R_0, q = W\cap N_0, a = W\cap N_7$. We have 
$3p\sim 2q+a$ and the fibration defines a $g_3^1$ on $W$ spanned by the divisors $3p$ and $2q+a$. 

We also observe that $X$ contains another elliptic fibration $|F'|$ with two reducible fibres 
$$F_1' = 3N_0+R_8+2S+3N_1+4N_2+5N_3+6N_4+4N_5+2N_6$$ and
$$F_2' = 2R_0+A+2R_1+3R_2+4R_3+3R_4+2R_5+R_6$$ 
of type $\tilde{E}_8$ and $\tilde{E}_8$. The curve $R_7$ is a section. The involution $\sigma$ switches the two fibres and induces the hyperelliptic involution on $W$. Its set of fixed points are 2 points on the curve $R_8$ and 6 points on $W$. Also note that $\sigma$ maps the fibration $|F|$ to the fibration $|F'|$.
It is easy to see that $q = \sigma(p)$. This gives $3p\sim 2K_W-2p+a$, hence
$$|5p| = |2K_W+a|.$$

\begin{teo}\label{cover} The linear system $|F'+F|$ defines a  map $f:X\to \bbP^1\times \bbP^1$ of degree 2. Its branch locus is a curve of bidegree $(4,4)$ which is equal to the union of a curve $B$ of bidegree $(3,3)$ and two rulings $E_1,E_2$. The curve $B$ has  2 cusps $q_1,q_2$. The   cuspidal tangent at $q_i$ is equal to $E_i$. The ramification curve of $f$ is equal to $W$. The automorphism $\tau$ of $X$ is the deck transformation of $f$, the quotient $X/(\tau)$ admits a birational morphism to $\bbP^1\times \bbP^1$ which resolves the singularities of the branch curve. The automorphism $\sigma$ is induced by the automorphism $\bar{\sigma}$ of $\bbP^1\times \bbP^1$ which switches the factors. It leaves the curve $B$ invariant and switches $E_1$ and $E_2$. 
\end{teo}

\begin{proof} We have $(F+F')^2 = 4$ and the restriction of $|F'+F|$ to a nonsingular fibre of each fibration is a degree 2 map. This easily implies that the linear system defines a degree 2 map $f$ to a nonsingular quadric in $\bbP^3$ and the pre-images of the rulings are the fibrations $|F|$ and $F'|$. The map $f$ blows down the curves
$R_0,R_1,R_2,R_3,R_4,R_6,N_0,N_1, N_2,N_3,N_4,N_5, N_6$. Its restriction  to $W$ is a birational map  defined by a 3-dimensional linear subsystem of
$|3p+3q| = |3K_W|.$ The rest of the assertions are easy to verify.
\end{proof}

Consider the automorphism $\bar{\sigma}$ of $\bbP^1\times \bbP^1$ and let 
$$\pi:\bbP^1\times \bbP^1\to \bbP^1\times \bbP^1/(\bar{\sigma}) \cong \bbP^1$$ be 
the natural projection to the orbit space.  Its locus of fixed points is the diagonal $\Delta$. The image of $\Delta$ on $\bbP^2$ is a conic $Q$. The image of $W$ is a cuspidal cubic $G$. The images of $E_1$ and $E_2$ is the cuspidal tangent $T$ of $G$. It is also a tangent of the conic  $Q$. The curves $Q$ and $G$ intersect at 6 points, the ramification points of the hyperelliptic involution of $W$.   

\begin{teo} Let $\bar{f}:Y\to \bar{Y}\to \bbP^2$ be a minimal resolution of the double cover $\bar{Y}$ of $\bbP^2$ branched along  
the union of the curves $G,Q$ and $T$. Then $Y$ is a Kummer surface birationally isomorphic to the quotient of $X$ by $\sigma$. 
\end{teo}

\begin{proof} We have a commutative diagram
\[\begin{CD}
X@>f>>\bbP^1\times \bbP^1\\
@V/\sigma VV@V/\bar{\sigma}VV\\
Y@>\bar{f}>>\bbP^2
\end{CD}
\]

It follows from \cite{Na} that $Y$ is a Kummer surface of $\text{Jac}(C)$ for some curve $C$ (not isomorphic to $W$).
\end{proof}

\begin{oss} One can reverse the construction of Theorem~\ref{cover}. Starting from a genus 2 curve $W$ together with a point $p\in W$ such that $|3K_W-p|\ne \emptyset$ and $p$ is not a Weierstrass point, we construct a 2-cuspidal model  $W$ in $\bbP^1\times \bbP^1$ as in the theorem. Then taking the double cover we get a K3 surface $X$ with Picard lattice containing $U\oplus E_8\oplus E_7$ and $U\oplus \langle -2\rangle \subset T_X$. Replacing $p$ by the conjugate point $q$ under the hyperelliptic involution, we get the same surface $X$. As was explained to me by J. Harris, the number of pairs $(p,q)$ as above is equal to 16. Thus we obtain that the moduli space of K3 surfaces marked with the lattice $U\oplus E_8\oplus E_7$ is isomorphic to a $(16:1)$-cover $\mathcal{M}_2'$ of the moduli space $\mathcal{M}_2$ of genus 2 curves, and on the other hand, via periods it is isomorphic to the moduli space $\mathcal{A}_2$ of principally polarized abelian surfaces. 
This defines a birational isomorphism between $\mathcal{M}_2'$ and $\mathcal{M}_2$.
\end{oss}

 \end{document}